\documentclass[11pt]{article}
\usepackage{amssymb}
\usepackage{amsthm}
\usepackage{amsmath}

\newtheorem{theorem}{Theorem}
\newtheorem{cor}[theorem]{Corollary}
\newtheorem{lemma}[theorem]{Lemma}
\newtheorem{prop}[theorem]{Proposition}

\newtheorem{claim}{Claim}
\newtheorem{definition}{Definition}
\newtheorem{question}{Open Question}

\newcommand{\restrict}[2]{#1\mspace{-2mu}\mathbin{\restriction}\mspace{-1mu} #2}

\newcommand{\F}{\mathbb{F}}

\newcommand{\N}{\mathbb{N}}
\newcommand{\Q}{\mathbb{Q}}

\newcommand{\Z}{\mathbb{Z}}

\newcommand{\CA}{\mathcal {A}}
\newcommand{\CB}{\mathcal {B}}

\newcommand{\CL}{\mathcal {L}}
\newcommand{\CM}{\mathcal {M}}
\newcommand{\CR}{\mathcal {R}}

\newcommand{\CV}{\mathcal {V}}

\newcommand{\age}{age}
\newcommand{\Aut}{Aut}
\newcommand{\Emb}{Emb}
\newcommand{\IGL}{IGL}

\newcommand{\osc}{\mathrm{osc}}
\newcommand{\supp}{\mathrm{supp}}
\newcommand{\maxsupp}{\mathrm{maxsupp}}
\newcommand{\minsupp}{\mathrm{minsupp}}
\newcommand{\FL}{\F^{[\lambda]}}
\newcommand{\FO}{\F^{[\omega]}}
\newcommand{\FqO}{\F_q^{[\omega]}}
\newcommand{\OO}{\overline{0}}

\newcommand{\blp}{\big(} 
\newcommand{\brp}{\big)} 

\begin{document}

\title{Partitions and Indivisibility Properties of Countable Dimensional Vector Spaces}

\author{C. Laflamme
\thanks{Supported by NSERC of Canada Grant
\# 690404}\\University of Calgary\\Department of
Mathematics and Statistics\\2500 University Dr. NW. Calgary
Alberta Canada T2N1N4\\
{\tt laf@math.ucalgary.ca}
\and
L. Nguyen Van Th\'e \thanks{The author would like to thank the support 
of the Department of Mathematics \& Statistics Postdoctoral Program at the University of Calgary}
\\Universit\'e de Neuch\^{a}tel
\\Institut de Math\'ematiques
\\Rue Emile-Argand 11
\\2007 Neuch\^{a}tel
\\Suisse 
\\ {\tt  lionel.nguyen@unine.ch}
\and
M.~Pouzet \\PCS, Universit\'e Claude-Bernard Lyon1,
\\ Domaine de Gerland -b\^at. Recherche [B], 50 avenue Tony-Garnier 
\\F$69365$ Lyon cedex 07, France
\\{\tt pouzet@univ-lyon1.fr }
\and
N.~Sauer\thanks{Supported by NSERC of Canada Grant \# 691325} 
\\Department of Mathematics and Statistics 
\\The University of Calgary, Calgary
\\Alberta, Canada T2N1N4
\\{\tt nsauer@math.ucalgary.ca}
}

\maketitle
\date{ }

\medskip

\begin{abstract} 
We investigate infinite versions of vector and affine space partition 
results, and thus obtain examples and a counterexample
for a partition problem for relational structures. In particular we
provide two (related) examples of an age indivisible relational structure which is
not weakly indivisible.
\end{abstract}



\section{Introduction}

In the present paper we study the divisibility properties of some
classical groups by studying the divisibility properties of infinite
dimensional vector spaces, generalizing those obtained for finite
dimensional spaces. Before we can state the main results, we review
basic divisibility notions from structural Ramsey theory.

A relational structure consists of a nonempty base set along with a
collection of finitary relations which are defined on it.  Let $\CR$
be a relational structure with base set $R$.  An induced substructure
of $\CR$ is a relational structure on a subset $\CR'$ of $\CR$
obtained by restricting all of the relations to $\CR'$. The {\em
skeleton} of $\CR$ is the set of finite induced substructures of
$\CR$, namely all the induced substructures of $\CR$ with base sets
nonempty finite subsets of $\CR$.  A typical example of a relational
structure is a graph consisting of vertices (the base set) together
with edges (one binary relation). Its induced substructures are what
are normally called (induced) subgraphs.
The {\em age} of $\CR$, $\age(\CR)$, is the class of finite
relational structures isomorphic to an element of the skeleton. A {\em
local isomorphism} of $\CR$ is an isomorphism between two elements of
the skeleton. We denote by $\Aut(\CR)$ the automorphism group of
$\CR$. The relational structure $\CR$ is {\em homogeneous} if every
local isomorphism of $\CR$ has an extension to an automorphism of
$\CR$. An {\em embedding} of $\CR$ into $\CR$ is an isomorphism of
$\CR$ to an induced substructure of $\CR$. We denote by $\Emb(\CR)$
the set of embeddings of $\CR$ into $\CR$.  For a set $A$ (possibly
not contained in $R$) we denote by $\restrict{\CR}{A}$ the relational
structure induced by $R \cap A$. 


We recall the following three notions from structural Ramsey theory
(see for example the Appendix of \cite{Fra}), together with the notion
of a {\em uniform partition} which we believe to be new (at least in
this form).

\begin{definition}

\begin{enumerate}

\item A relational structure
$\CR$ is {\em indivisible} if for every partition $(P_0, P_1)$
of $R$ there exists an element $\epsilon\in \Emb(\CR)$ and an
$i\in 2$ with $\epsilon[R]\subseteq P_i$.

\item A relational structure $\CR$ is said to have a {\em uniform
partition} if there is a finite partition $(U_i: i  \in  n)$ of $R$ such
that 
\begin{enumerate}
\item for all $\epsilon\in \Emb(\CR)$ and all $i \in n$,  $\age(\restrict{\epsilon[\CR]}{U_i})
= \age(\CR)$, 
\item for every partition $(P_0, P_1)$ of $R$ and
any $i \in n$, there exists an element $\epsilon\in \Emb(\CR)$ and
$j \in 2$ with $(\epsilon[R] \cap U_i) \subseteq P_j$.
\end{enumerate}

\item  The relational structure $\CR$ is {\em weakly indivisible} if
for every partition $(P_0,P_1)$ of $R$ with
$\age(\restrict{\CR}{P_0})\neq
\age(\CR)$ there exists an element $\epsilon\in
\Emb(\CR)$ with $\epsilon[R]\subseteq P_1$. 

\item The relational
structure $\CR$ is {\em age indivisible} if for every partition
$(P_0,P_1)$ of $R$, there is an $i \in 2$ with
$\age(\restrict{\CR}{P_i})= \age(\CR)$.  
\end{enumerate}
\end{definition}

\noindent Each of the above properties easily implies the one following it. 
The more familiar notion of a canonical partition is defined as in
2 above by simply replacing the first condition by the weaker
requirement that any copy of $\CR$ meets every block of the partition (see \cite{Sa2}):

\begin{definition}
A relational structure $\CR$ is said to have a {\em canonical 
partition} if there is a finite partition
$(U_i: i  \in  n)$ of $R$ such that 
\begin{enumerate}
\item  for all $\epsilon\in
\Emb(\CR)$ and all $i \in n$,  
$\epsilon[\CR] \cap U_i \not= \emptyset$ 
\item for every partition
$(P_0, P_1)$ of $R$ and any $i \in n$, there exists an element
$\epsilon\in \Emb(\CR)$ and $j \in 2$ with $(\epsilon[R] \cap U_i)
\subseteq P_j$.
\end{enumerate}
\end{definition}

A uniform partition is clearly a canonical partition, and conversely a
canonical partition together with the weak indivisibility property
implies that it is a uniform partition. However the existence of a
canonical partition alone does not suffice; indeed there are countable
homogeneous structures having a canonical partition but not age
indivisible, namely the structure consisting of two disjoint copies of
the rationals equipped with the usual linear order on each copy (and
no other relation).

More generally it is known that most of the above implications are not 
reversible. Indeed there are countable homogeneous divisible (meaning
not indivisible) structures having uniform partition; these will be
further described in detail in a forthcoming paper (see
\cite{Sa3}). There are examples (see \cite{Sa1}) of countable weakly
indivisible homogeneous divisible structures. We will see an example
of a homogeneous structure below which is weakly indivisible but does
not have a canonical partition, and therefore not a uniform one, and
therefore is divisible.

We were not aware of any example of a countable age indivisible
homogeneous structure which is not weakly indivisible. See \cite{Sa3}
for the fact that the standard age indivisible examples of countable
homogeneous relational structures are also weakly
indivisible. Actually, there was no example of an age indivisible
structure (homogeneous or not) which is not weakly indivisible. We
provide two examples below obtained from vector spaces. They are
unfortunately not as simple as one would intend in the sense that one
has infinitely many relations. Although the other has a single
relation, it is not homogeneous. It remains open whether it is possible to produce
such an example of a countable homogeneous structure having only
finitely many relations.

\medskip

One of the standard tools to prove age indivisibility of relational
structures is the Hales-Jewett Theorem \cite{HJ}; weak indivisibility
seems then related to an infinite version of that theorem at least in
cases similar to the vector space situation discussed in the present
paper, and we shall therefore be interested in infinite dimensional
vector spaces. In the case of finite dimensional vector spaces over a
finite field $\F_q$, Graham, Leeb and ~Rothschild, proved the
following.

\begin{theorem}\label{theorem:glr}\cite{GLR}
For all $d,k,t \geq 0$, there exists $n=GLR^t(d,k)$ with the
property that for any $n$-dimensional vector space $V$ over $\F_q$
and any  colouring of all $t$-dimensional (affine) subspaces into
$k$ colours, there exists a $d$-dimensional (affine)  subspace $U \subset
V$ such that all its $t$-dimensional (affine) subspaces have the same colour. 
\end{theorem}

\noindent  The reason for writing the adjective ``affine'' in parenthesis 
is that the above result where the notion of subspaces is interpreted
as ``affine'' is equivalent, as proved by Graham and Rothschild, to
the corresponding one using the normal sense of subspace. In that
latter sense of usual vector subspace, the Theorem for $t=0$ has no
content. But for $t=1$ (or equivalently its affine version for $t=0$)
is already powerful , it implies in particular the following particular
case which will be used in the proof of Theorem
\ref{theorem:infiniteblueaffine}:

\begin{cor} \label{cor:lines}[t=1]
For all $d,k \geq 0 $, there exists $n=GLR(d,k)$ such that for any
$n$-dimensional vector space $V$ over $\F_q$ and any colouring of the
lines of $V$ into $k$ colours, there exists a $d$-dimensional subspace
$W \subset V$ all of whose lines have the same colour.
\end{cor}

In the affine case we get the following. 

\begin{cor} \label{cor:vectors} [t=0 (affine)]
For all $d,k \geq 0 $, there exists $n$ such that for any
$n$-dimensional vector space $V$ over $\F_q$ and any colouring of $V$
into $k$ colours, there exists a monochromatic $d$-dimensional affine
subspace $W \subset V$.

\end{cor}

A central motivation for this research is to investigate infinite
versions of this result. To do so, we shall interpret a vector space
as a relational structure, in which the Ramsey partition properties
described above correspond to affine version of the usual vector
subspace and related notions, and we assume this notion for the
remainder of the paper unless specifically mentioned otherwise.  In
particular we will be interested in the affine transformations of $V$
onto $V$, forming a group called the (Inhomogeneous) General Linear
group of $V$, and denoted by $\IGL(V)$. 

\medskip 

Hindman showed in \cite{hindman} that a vector space of countable
dimension over $\F_2$ is indivisible. On the other hand we provide a
proof in Theorem~\ref{theorem:nouniformpartition} of the well known fact
that a vector space of countable dimension over any other field is
divisible, in fact does not have a uniform partition. We also prove in
Theorem~\ref{theorem:infiniteblueaffine} that over any finite field, a
vector space of countable dimension is weakly indivisible.  Over an
infinite field, we shall show that a vector space of countable
dimension is not weakly indivisible, in fact it can be divided into
two parts such that none of the parts contains an affine line (see
Theorem~\ref{theorem:infi}). On the other hand, all infinite
dimensional vector spaces are age indivisible. So a countable
dimensional vector space over the rationals provides an example
of an age indivisible, not weakly indivisible relational structure
(with infinitely many relations).

It is known and we will provide a proof in
Lemma~\ref{lem:affineembedding} that if $V$ is a vector space over the
rationals $\mathbb{Q}$ and $\epsilon: V\to V$ an injection with
$\epsilon(\frac{a+b}{2}) = \frac{\epsilon(a) +
\epsilon(b)}{2}$ for all $a,b\in V$ then $\epsilon\in \Emb(V)$, that
is $\epsilon$ is an affine transformation.  It follows that if $\CM_V$
is the relational structure with base set $V$ and ternary relation
$\mu(a,b,c)$ if and only if $b=\frac{a+c}{2}$ then
$\Emb(V)=\Emb(\CM_V)$ and in particular that $\Aut(\CM_V)=\IGL(V)$. We
prove in Theorem~\ref{theorem:notweakly} that $\CM_V$ does constitute
another example of an age indivisible but not weakly indivisible
relational structure.  Unfortunately $\CM_V$ is not homogeneous: for
take any $n$ element sequence on an affine line no three of whose
points are in the midpoint relation $R$. Then there is a local
isomorphism $\alpha$ to any other such $n$ element sequence on any
other line, even if $\alpha$ is not an affine transformation, in which
case $\alpha$ cannot be extended to an element of $\IGL(V)$.

\medskip 

It is worth noting that the above Ramsey properties for homogeneous
structures are properties of the automorphism group seen as a
permutation group. Conversely, given a permutation group closed in the
product  topology, there exist homogeneous structures with the given
group as automorphism group (see \cite{Fra} for some discussion on
this). That is, those divisibility properties can be studied as
permutation group properties.

\section{The Affine Space Structure} \label{section:vsstructure} 

An affine transformation $\alpha$ of a vector space $V$ onto another
vector space $W$ is a function of $V$ to $W$ for which there exists an
element $w\in W$ and an invertible linear transformation $\rho:
V\to W$ so that $\alpha(v)=w+\rho(v)$ for all $v\in V$. The affine
transformations of $V$ to $V$ form a group $\IGL(V)$, the
Inhomogeneous General Linear group of $V$, and we denote the set of
corresponding affine embeddings of $V$ into $V$ by $\Emb(V)$.

\noindent A sequence $\langle v_i : i \in n \rangle$ of affinely dependent
elements of $V$ is called an {\em affine cycle} of $V$ if no proper
subsequence is affinely dependent. Two affine cycles $\langle v_i : i \in n \rangle$
and $\langle v^\prime_i : i \in n \rangle $ of $V$ are {\em equivalent} if there is
an invertible affine transformation $\tau$ of the affine space
generated by $\{v_i: i\in n\}$ to the affine space generated by
$\{v_i^\prime: i\in n\}$ with $\tau(v_i)=v^\prime_i$ for all $i\in
n$. Interpreting every equivalence class of affine cycles of $V$ as a
relation on $V$ yields a homogeneous relational structure $\CV$.
Under this situation it turns out that $\Aut(\CV)=\IGL(V)$, and that
the set of affine embeddings of $V$ is equal to the set of embeddings
of $\CV$. The relational structure $\CV$ is what we call the {\em
affine cycle structure of $V$}.  Every element of the skeleton of
$\CV$ generates affinely a finite dimensional affine subspace of $V$
and hence every element of the age of $\CV$ can be affinely embedded
into a finite dimensional affine subspace of $V$. Of course when the
field is finite, then every finite dimensional affine subspace of $V$
is also part of the skeleton.  As a consequence we obtain the
following translation of the structural Ramsey properties described
above applied to $\CV$. Thus we shall say by abuse of terminology that
a vector space $V$ is {\em indivisible} if for every partition
$(P_0,P_1)$ of $V$ there exists an affine embedding $\epsilon$ of $V$
and an $i\in 2$ such that $\epsilon[V]\subseteq P_i$. The vector space
$V$ has a {\em uniform partition} if there is a finite partition
$(U_i: i \in n)$ of $V$ such that (1) for all affine embeddings
$\epsilon$ of $V$ and all $i \in n$,
$\age(\restrict{\epsilon[V]}{U_i}) = \age(V)$, (2) for every partition
$(P_0, P_1)$ of $V$ and any $i \in n$, there exists an affine
embedding $\epsilon$ of $V$ and $j \in 2$ with $(\epsilon[V] \cap U_i)
\subseteq P_j$. The vector space $V$ is {\em weakly indivisible} if
for every partition $(P_0,P_1)$ of $V$ with $\age(V) \neq
\age(\restrict{V}{P_0})$, there exists an affine embedding $\epsilon$
of $V$ with $\epsilon[V]\subseteq P_1$.  A vector space $V$ is {\em
age indivisible} if for every partition $(P_0,P_1)$ of $V$, there is
an $i \in 2$ with $\age(\restrict{V}{P_i}) =
\age(V)$.

Using a standard  compactness argument one can show that a relational
structure $\CR$ is age indivisible if and only if the age of $\CR$ is
a Ramsey family. That is, if for every element $\CA$ in the age of
$\CR$ with base $A$ there exists an element $\CB \in \age(\CR)$ with
base set $B \supseteq A$ so that for every partition $(P_0,P_1)$ of
$B$ there exists an embedding $\epsilon$ and an $i\in 2$ with
$\epsilon[A]\subseteq P_i$. It follows readily from this that two
relational structures with the same age are either both age
indivisible or neither is. This characterization is also useful in
proving age indivisibility; indeed if the age of a relational
structure $\CR$ is closed under products (for an appropriate
definition of product) as for vector spaces in the case of this paper,
then it follows from the Hales-Jewett Theorem (see
\cite{HJ}), that the age of $\CR$ has the Ramsey property and
hence that $\CR$ is age indivisible. In the case of vector spaces over
a finite field $\F_q$, one can use $\F_q$ itself as the required
alphabet and choose a sufficiently large product and observe that a
combinatorial line is an affine line, and then more
generally that a combinatorial space is an affine
space. It follows from this discussion that the affine cycle structure
of a vector space is age indivisible.


\medskip

We conclude this section by reviewing some basic notions and notation
for vector spaces that will be used throughout this paper. If
$\lambda$ denotes the dimension of a vector space $V$ over a field
$\F$, then we identify $V$ with $\FL$ as the set of functions
$h:\lambda \rightarrow \F$ taking nonzero values only finitely many
times.  The \emph{support} of such a function $h$ is the set
$\supp(h):=\{\alpha \in \lambda: h(\alpha)\not =0\}$.  We denote by
$\OO\in V$ the constant sequence with value $0\in \F$, with the
understanding that $\supp(\OO)=\emptyset$.  We write $\maxsupp(h)$ for
the largest element of $\supp(h)$ if $h\not =\OO$, and
$\maxsupp(\OO):=-\infty$; finally, we set $\hat{h}:=h\blp \maxsupp(h)\brp$ if
$h\not=\OO$ and $\hat{\OO}:=0$. Similarly we write $\minsupp(h)$
for the smallest element of $\supp(h)$ if $h\not =\OO$, and
$\minsupp(\OO):=-\infty$, and set $\check{h}:=h\blp \minsupp(h) \brp$ if
$h\not=\OO$ and $\check{\OO}:=0$. More generally the support
$\supp(A)$ of a subset $A$ of $\FL$ is the set $\bigcup_{h\in
A}\supp(h)$.

For two finite subsets $X$ and $Y$ of $\lambda$, we write $X\lll Y$ if
the maximum of $X$ is strictly smaller than the minimum of $Y$. We extend this
notation to $f,g\in \FL$ by writing $f\lll g$ if $\supp(f)\lll\supp(g)$,
and to $A\lll B$ for subsets $A$ and $B$ of $\FL$ if $f\lll g$ for all
$f\in A$ and $g\in B$.

A subset ${A}$ of $V$ is an {\em affine subspace} if there is a
(unique) subspace $W$ and an element $v$ of $V$ with
$A=v+W:=\{v+w: w\in W\}$. The dimension of the affine subspace $A$ is
the dimension of $W$.  A (affine) {\em line} is a (affine) subspace of
$V$ of dimension one, and we denote by $\CL$ the set of lines of $V$.
Every line of $V$ contains exactly one element $f$ with $\hat{f}=1$,
and conversely every $f \in V$ with $\hat{f} =1$ generates a line
$\{af \mid a\in \F \}$ which we denote by $\langle f \rangle$;
that is we name a line by its unique element $f$ with $\hat{f}=1$.


\medskip 

We shall be mostly interested in the countable case $\lambda=\omega$,
but many of the results presented generalize to vector spaces of an
arbitrary dimension $\lambda$.

\section{Vector spaces of countable dimension}

For this section, fix a vector space $V$ of countable dimension over an arbitrary 
field $\F$, and as described above we may assume that $V=\FO$. We
begin  by producing a manageable and interesting structure for an
infinite dimensional (affine) subspace.

\begin{lemma}\label{lem:xnotinsupp}
Let $W$ be an infinite dimensional subspace of $V$ and $x \in
\omega$.  Then there  exists an infinite  dimensional subspace $U$  of
$W$ so that $x \lll \supp(f)$ for all $f\in U$.
\end{lemma}

\begin{proof}
By a repeated application, it suffices to prove that for each $x \in
\omega$, there is a nonzero vector  $f \in W$ with $x \lll \supp(f)$.

\noindent Since $W$ is infinite dimensional, there must be two linearly 
independent vectors $f,g \in W$ such that $f \restriction x$ and $g
\restriction x$ are linearly dependent. That is $a f \restriction x +b g
\restriction x = \OO$ for some $a,b \in \F$ not both zero. But then $x \lll  af + bg \neq \OO$. 
\end{proof}

An iterated application of Lemma~\ref{lem:xnotinsupp}
allows to construct the following  structure for an
infinite dimensional affine subspace of $V$. 

\begin{prop}\label{prop:infiniteblocks}
Every infinite dimensional affine subspace $v+W$ of $V$ contains an
infinite sequence $(v+f_i: i\in \omega)$ such that 
\[ v \lll f_i\lll f_{i+1} \mbox{  for all } i \in \omega .\]
\end{prop}

\subsection{Indivisibility}


In \cite{hindman}, Hindman proved that a vector space of countable
dimension over $\F_2$ is indivisible. This fact is an immediate
consequence of  (and equivalent to) Hindman's finite union Partition Theorem using
Proposition~\ref{prop:infiniteblocks}.

\begin{theorem} \cite{hindman} 
If $V$ is a vector space over $\F_2$ of countable dimension, then
$V$ is indivisible. 
\end{theorem}


\subsection{Uniform Partitions}

It is a well known folklore result that every countable dimensional
vector space over any other field than  $\F_2$ is divisible, in fact does not contain a
uniform or even a canonical partition.

\begin{theorem} (Folklore) \label{theorem:nouniformpartition}
If $\F\not=\F_2$, then any countable dimensional vector space $V$ over
$\F$ is divisible. 
\noindent In fact $V$ does not have a uniform or  even a canonical partition.
\end{theorem}

\begin{proof}
For $f \in V=\FO$, define $\osc(f)$ as the number of times that $f$
changes from a nonzero value to a different nonzero value as we cover
the support of $f$. That is, if $\supp(f)=\{x_i: i \in n\}$ is listed in
increasing order, then
\[ \osc(f) = | \{i \in n-1: f(x_i) \not= f(x_{i+1}) \}|. \]
Observe that if $f \lll g$, then $\osc(f+g)\geq \osc(f)+\osc(g)$, with
equality iff the last value of $f$ equals the first value of $g$ (if
those values are different, $\osc(f+g)= \osc(f)+\osc(g)+1$). 

\noindent Now consider a sequence $\langle f_i : i  \in  n \rangle$ from a subspace $W$ 
of $V$ such that $f_i \lll f_{i+1}$. By an appropriate scalar
multiplication, we may assume that $\hat{f_i} = \check{f}_{i+1}$ for
all $i \in n-1$, and therefore $s := \osc \left( \sum_{i \in n} f_i
\right) = \sum_{i \in n} \osc(f_i)$. Since $\F\not=\F_2$, 
choose for each $i \in n-1$ an $a_i \not= 0 \in \F$ such that $a_i
\not= a_{i+1}$. But now observe that for each $j \in n$,
\[ \osc\left( \sum_{i=0}^{j} a_i f_i  + a_{j+1} \sum_{i=j+1}^{n-1} f_i\right)= s+j+1 .\]
This means that on any infinite dimensional subspace or even affine
subspace of $\FO$, the range of the oscillation function contains
arbitrarily long intervals. Hence, there cannot be any canonical
partition.
\end{proof}








\subsection{Weak Indivisibility}

\subsubsection{Weak Indivisibility over Finite Fields}

Let $\mathbb{F}_q$ be the finite field of $q$ elements. We shall prove
that $V_q = \FqO$ is weakly indivisible.

\begin{lemma}\label{lem:xnotinsupp2}
Let $k\in \omega$, $W$ a subspace of $V_q$ of dimension at least $k+1$, and $x
\in \omega$ arbitrary. Then there exists a $k$-dimensional subspace
$U$ of $W$ such that $x \notin\supp(U)$.
\end{lemma}

\begin{proof}
The map $W \rightarrow \mathbb{F}_q$ given by $f \rightarrow f(x)$ is
a linear map from a subspace of dimension at least $k+1$ into a
1-dimensional one, so its kernel must have dimension at least $k$. 
\end{proof}

\begin{cor}\label{cor:valuei}
Let $k \in \omega$, $W$ be a subspace  of $V_q$ of dimension at least $k+1$, and
$x \in\supp(W)$. Then for any $a\in \mathbb{F}_q$, there exists an affine 
$k$-dimensional subspace $A$ of $W$ so that $f(x)=a$ for all $f\in A$.
\end{cor}

\begin{proof}
Let $w \in W$ such that $w(x) \not= 0$, and by Lemma~\ref{lem:xnotinsupp2}
let $U$ a $k$-dimensional subspace of $W$ such that $f(x)=0$ for all
$f\in U$.

\noindent Then $A= \{a(w(x))^{-1}w + f : f \in U \}$ 
is the desired affine $k$-dimensional subspace.
\end{proof}

A colouring  of a subset $W$ of $V_q$ is called {\em end-determined } 
if, for every $a \in \mathbb{F}_q\setminus 0$, the set $\{f\in W :
\hat{f} =a \}$ is monochromatic.

\begin{lemma}\label{lem:endext}
Let $k \in \omega$, $W$ be a subspace of $V_q$  of dimension at least
$k+1$, and $v\in V_q$ with $\supp(v)\lll \supp(W)$.  Let $\Delta$ be a
colouring of the affine space $v+W$ into $red$ and $blue$ elements so
that $\Delta(v)=blue$, and so that $v+W$ does not contain a monochrome
$red$ affine subspace of dimension $k$.

\noindent If the colouring $\Delta$ is end-determined on $v+W$, 
then every element of $v+W$ is blue.
\end{lemma}

\begin{proof}
Assume for a contradiction that $h\in v+ W$ is red, and let
$a=\hat{h}$.  Then every $g\in v+ W$ with $\hat{g}=a$ is red since
$\Delta$ is end-determined on that space. Hence may assume without
loss of generality that $\maxsupp(h)=\maxsupp(W)$.

\noindent According to Corollary~\ref{cor:valuei}, there exists an affine
$k$-dimensional subspace $A$ of $W$ so that $\hat{f}=a$ for every
$f\in A$. But then $v+A$ is an affine red subspace of dimension $k$, a
contradiction.
\end{proof}

\begin{lemma}\label{lem:bluebasic}
Let $k \leq d  \in \omega$, $v\in V_q$, and $V$ a
$GLR(d+1,2^{q-1})$-dimensional subspace of $V_q$ with $v\lll V$.  Let
$\Delta$ be a colouring of the affine space $v+V$ into $red$ and
$blue$ elements so that $\Delta(v)=blue$ and $v+V$ does not contain a
monochrome $red$ affine subspace of dimension $k$. 

\noindent Then there exists a $d$-dimensional subspace $U$ of $V$ so that every
element of the affine space $v+U$ is blue.
\end{lemma}

\begin{proof}
Colour every line $L=\langle f \rangle=\{af : a\in \mathbb{F}_q\}$ (where $\hat{f}=1$) of
$V$ with the function $\gamma_L: \mathbb{F}_q\setminus\OO\to
\{red,blue\}$ given by $\gamma_L(a)=\Delta(v+af)$;  that
is with one of $2^{q-1}$ possible colours. Denote by $\Gamma$ this
colouring of the set of lines in $V$ with $2^{q-1}$ colours.

\noindent Then by Theorem~\ref{theorem:glr} there  exists  a
$d+1$-dimensional subspace $W$ of $V$ and a function $\gamma:
\mathbb{F}_q\setminus\OO\to \{red,blue\}$ so that
$\Gamma(\langle f \rangle )=\gamma$ for every line $\langle f \rangle \in
W$. But this means that the colouring $\Delta$ is end-determined on the
affine subspace $v+W$, which is therefore by assumption and
Lemma~\ref{lem:endext} monochrome blue.

\end{proof}

Toward the proof of our next result, we define recursively the number
$\Pi_n(d)$ for $n, d \in \omega$ by $\Pi_1(d):=GLR(d+1,2^{q-1})$ and
$\Pi_{n+1}(d):=GLR\big(\Pi_n(d)+1,2^{q-1}\big)$.

\begin{lemma}\label{lem:bluebasicext}
Let $k,d\in\omega$, $V$ be an infinite dimensional subspace of $V_q$,
and $v+A$ be an affine finite dimensional subspace of $V$.

\noindent Let $\Delta$ be a colouring of $V$ into $red$ and $blue$ elements so that
there is no monochrome red affine subspace of dimension $k$,  and so that
$\Delta$ is monochrome blue on $v+A$.

\noindent Then there exists a $d$-dimensional subspace 
$W$ of $V$ with $v+A\lll W$ and so that every element in $(v+A)+W$ is
blue.
\end{lemma}

\begin{proof}
List the elements of $v+A$ as $f_{n-1},f_{n-2}, f_{n-3},\dots, f_0$,
and using Lemma~\ref{lem:xnotinsupp}, let $W_{n-1}$ be a
$\Pi_n(d)$-dimensional subspace of $V$ with $v+A\lll W_{n-1}$.

Then $v+f_{n-1}+W_{n-1}$ is an affine space with $v+f_{n-1}$
blue. According to Lemma~\ref{lem:bluebasic} there there exists a
$\Pi_{n-1}(d)$-dimensional subspace $W_{n-2}$ of $W_{n-1}$ so that
$(v+f_{n-1})+W_{n-2}$ is monochrome blue.

More generally, assume that we have $1\leq i<n$ and a
$\Pi_{n-i}(d)$-dimensional subspace $W_{n-(i+1)}$ of $W_{n-i}$ so that
the space $(v+f_{n-i})+W_{n-(i+1)}$ is monochrome blue. Then
$v+f_{n-(i+1)}+W_{n-(i+1)}$ is an affine space with $v+f_{n-(i+1)}$
blue. According to Lemma~\ref{lem:bluebasic} there there exists a
$\Pi_{n-(i+1)}(d)$-dimensional subspace $W_{n-(i+2)}$ of $W_{n-(i+1)}$
so that the space $(v+f_{n-(i+1)})+W_{n-(i+2)}$ is monochrome blue.

\noindent We continue and for $i=n$ obtain a $d$-dimensional 
subspace $U$ of $W_{n-1}$ so that for every $v+f\in v+A$ and every
$g\in U$ the element $v+f+g$ is blue.

\end{proof}  

We now come to the main result of this section: $V_q$ is weakly
indivisible.

\begin{theorem}\label{theorem:infiniteblueaffine}
Let $V$ be a countable dimensional subspace of $V_q$, $k\in
\omega$, and $\Delta$ a colouring of $V$ into $red$ and $blue$ elements
so that $V$ contains no monochrome $red$ affine $k$-dimensional subspace.

\noindent Then there exists a monochrome blue affine 
subspace of $V$ of infinite dimension.
\end{theorem}

\begin{proof}
The space $V$ must contain at least one blue element $v$. Then $\{v\}$
is a 0-dimensional subspace which is monochrome blue. We obtain the
blue affine subspace of infinite dimension by repeated applications of
Lemma~\ref{lem:bluebasicext}.
\end{proof}


\subsubsection{Weak Indivisibility over Infinite Fields}

If the field $\F$ is infinite,
Theorem~\ref{theorem:nouniformpartition} has the following
strengthening, namely that $\FO$ is not weakly indivisible.

\begin{theorem}\label{theorem:infi}
Every countable dimensional vector space $V$ over an infinite field
$\F$ is not weakly indivisible. 

\noindent  In fact $V$ can be divided into two parts so that neither part 
contains an affine line. 
\end{theorem}

\begin{proof}
If $\F$ is a field of infinite size $\kappa$, then since the space is
of countable dimension we can enumerate the affine lines of $V$ as
$\langle L_\alpha : \alpha \in \kappa \rangle$. The intended set $A$ will be
constructed recursively as a sequence $\langle a_\alpha : \alpha \in
\kappa \rangle  \subseteq V$ such that for every $\alpha\in
\kappa$, $a_\alpha \in L_\alpha$ and such that no affine line intersects
$A_\alpha:=\{a_\beta: \beta \in \alpha\}$ in more than two points.

\noindent To do so, we pick $a_0 \in L_0$ arbitrary. Having defined $A_\alpha$, 
let $\CL_\alpha$ be the set of affine lines containing two distinct
points of $A_\alpha$, if any. If $L_\alpha$ already intersects
$A_\alpha$, let $a_\alpha \in (L_\alpha \cap A_\alpha)$. If not,
observe that since any two distinct affine lines intersect in at most
one point, $\CL_\alpha$ must have size less than $\kappa$, and
therefore we can choose $a_\alpha \in L_\alpha \setminus (\bigcup
{\CL_\alpha})$.

\end{proof}

In \cite{baumgartner}, Baumgartner proved the analog result for a
vector space of any dimension over the field of rational numbers.

\medskip 

Thus any countable dimensional  vector space over an infinite field
provides an example of an age indivisible but not weakly indivisible
homogeneous relational structure.



\section{Midpoint Structure}

In this section, we shall produce a somewhat simpler example of a
countable age indivisible but not weakly indivisible relational
structure. The structure will have a single ternary relation, but is
not homogeneous.

Before we proceed, let $M$ be a commutative monoid. An
\emph{arithmetic progression of length $n$} in $M$ is a sequence of
the form $(a+ix)_{i\in n}$ for some $a\in M$ and $ x\in M\setminus
\{0\}$. An \emph{infinite arithmetic progression} is defined
similarly. Clearly, an arithmetic progression of length $3$ is a
sequence of three elements $a_0, a_1, a_2$ of $M$ such that
$2a_1=a_0+a_2 $. Set $\mu_{M}:=\{(x,y, z)\in M^3: 2y=x+z\}$ and let
$\CM_M:= (M, \mu_M)$. In the case $M=\N$, we make the convention that
$\mu_{\N}$ denotes the set of triples associated to the additive
monoid on the nonnegative integers.

\noindent In the case where $M=V$ is a vector space over $\Q$, which will be the
main case of interest, then $\mu_{V}$ denotes the set of triples
associated with the addition on $V$.  Notice that in this case an
arithmetic progression of length $3$ is a sequence of three elements
$a_0, a_1, a_2$ where $a_1$ is the midpoint of the segment joining
$a_0$ and $a_2$. The ternary relational structure $\CM_V =
(V,\mu_V)$ is the {\em midpoint structure} associated with $V$. We
shall show in particular that $\CM_V$ is age indivisible but
not weakly indivisible. 

We first characterize the embeddings of such a structure
$\CM_V$ as simply the affine embedding of the underlying vector
space.

\begin{lemma} \label{lem:affineembedding}
Let $V$ and $V'$ be two vector spaces over $\Q$. A map
$\alpha:V\rightarrow V'$ is an embedding of\/ $\CM_V$ into
$\CM_{V'}$ if and only if it is an affine embedding of the
underlying vector spaces. 
\end{lemma}

\begin{proof}
First an affine embedding $\alpha$ of $V$ into $V'$ does satisfy
$\alpha(\frac{a+b}{2})=\frac{\alpha(a)+\alpha(b)}{2}$ for all $a,b\in
V$, and is therefore an embedding of $\CM_V$ into $\CM_{V'}$.  

Conversely this condition implies that $\alpha$ is an affine
transformation, indeed it suffices to show that $\beta: V\to V^\prime$
given by $\beta(a)=\alpha(a)-\alpha(\OO)$ for all $a\in V$ is a linear
transformation. 

\noindent Note that:
\[
\frac{\alpha(a)+\alpha(b)}{2}
	=\alpha(\frac{a+b}{2})=\alpha(\frac{a+b+\OO}{2})
	=\frac{\alpha(a+b)+\alpha(\OO)}{2}.
\]
Hence $\alpha(a+b)=\alpha(a)+\alpha(b)-\alpha(\OO)$, and therefore
$\beta(a+b)=\alpha(a+b)-\alpha(\OO)=\alpha(a)+\alpha(b)-2\alpha(\OO)=\beta(a)+\beta(b)$. It
follows immediately that $\beta(xa)=x\beta(a)$ for all rational
$x$. Therefore $\beta$ is a linear transformation, and since moreover
$\alpha$ is one to one, then it is an affine embedding as desired.
\end{proof}

We therefore immediately have the following Corollary.

\begin{cor}
The group $Aut(\CM_V)$ of automorphisms of $\CM_V$ and
the Inhomogeneous General Linear group $IGL(V)$ of $V$
coincide. Moreover, the closure of $Aut(\CM_V)$ (in teh product topology) consists of
exactly the affine embeddings of $V$ into itself, namely the set
$\Emb(V)$.
\end{cor}

The age indivisibility of $\CM_V$ will follow from the
following Lemma.

\begin{lemma}\label{lem:age} 
If $G$ is a torsion free abelian group then $\CM_G$ and
$\CM_{\N}$ have the same age.
\end{lemma}

\begin{proof} 
Since $\N$ can be identified with a subgroup of $G$, the age of
$\CM_{\N}$ is a subset of the age of $\CM_G$.

Conversely let $F$ be a finite subset of $G$. We will show that there
is a map $f: F\rightarrow \N$ which is an embedding of
$\restrict {\CM_G} {F}$ into $\CM_{\N}$. 

\noindent  Since $G$ is torsion free, the subgroup  generated by $F$ is 
isomorphic to a finite direct sum of the integers $\Z$. Without loss of generality,
we may therefore suppose that $F\subseteq \Z_+^n$, and let $k \in \N$
be large enough such that $F \subseteq k^n$. Finally define $\sigma:
k^n \rightarrow \N$ by
\[ \sigma( \langle x_i: i \in k \rangle ) = \sum_i x_i k^i . \]
Then $\sigma \restriction F$ is easily seen to be the required
embedding.
\end{proof}

The following tool is key in showing the failure of weak
indivisibility, and the proof is similar to that of Theorem \ref{theorem:infi}.

\begin{lemma}\label{arithmetic}
A countable commutative monoid $M$ in which equations of the
form $a+x=b$ and $a+2x=b$ have only finitely many solutions in $x$ for all $a,b\in
M$ contains a subset having no three elements
forming an arithmetic progression, but containing a point from  every
infinite arithmetic progression.
\end{lemma}

\begin{proof}
Note that the set of infinite arithmetic progressions is also
countable. Also note that if $S\subseteq M$ is finite and $X$ is an
infinite arithmetic progression with $X\cap S=\emptyset$ then there is
an element $x\in X\setminus S$ which does not form a three element
arithmetic progression with any two of the elements in $S$. This
follows from the fact that, for every two elements $\{a,b\}$ in $S$,
there are only finitely many equations of the above type in $M$
potentially producing a three element arithmetic progression with
$\{a,b\}$, each with only finitely many solutions by
hypothesis. Namely for any elements $x,y \in M$ with $a+x=b$ or
$b+y=a$ respectively, each of $a+2x$ and $b+2y$ will form an
arithmetic progression with $\{a,b\}$; moreover for any such $x,y$ as
above, for any elements $x',y' \in M$ with $x'+x=a$ or $y'+y=b$
respectively, $x'$ and $y'$ will also form an arithmetic
progression with $\{a,b\}$; finally for any elements $x,y \in M$ with
$a+2x=b$ or $b+2y=a$ respectively, $a+x$ and $b+y$ will form an
arithmetic progression with $\{a,b\}$.

Now enumerate the elements of $M$ into the $\omega$-sequence $x_0,
\dots, x_n, \dots$ and the set of infinite arithmetic progressions
into the $\omega$-sequence $X_0, \dots, X_n, \dots$. We construct the
sequence $y_0, \dots,y_n, \dots $ such that for every integer $n$, the
set $Y_n:=\{y_i: i\in n\}$ contains no three elements forming an
arithmetic progression, but meets $X_i$ for every $i \in n$. The element
$y_0$ is an arbitrary element in $X_0$. If $Y_n$ is already
constructed and $X_n\cap Y_n\not=\emptyset$, then let $y_n\in X_n\cap
Y_n$. If on the other hand $X_n\cap Y_n=\emptyset$, then let $y_n\in
X_n$ such that it does not form a three element arithmetic progression
with any pair of elements in $Y_n$. This completes the proof.
\end{proof}

We are now ready to prove the main result of this section.

\begin{theorem}\label{theorem:notweakly} 
Let $V$ be a vector space of countable dimension over $\Q$ and
$\CM_V:= (V, \mu_V)$ be the midpoint structure
associated with the vector space $V$. Then:
\begin{enumerate}
\item $\CM_V$ is age indivisible.
\item $\CM_V$ is not weakly indivisible.
\item $\CM_V$ is universal for its age: every countable 
$\CR:=(R, \mu')$ with the same age as $\CM_V$ is embeddable
into $\CM_V$.
\end{enumerate}

\end{theorem}
\begin{proof}
We prove each part separately. 

\noindent Item 1. We already observed in Section \ref{section:vsstructure} 
that two relational structures with the same age either are both age
indivisible or both age divisible.

\noindent According to Lemma~\ref{lem:age}, $\CM_V$ and $M_{\N}$ have the
same age. Thus it suffices to prove that for each subset $A$ of $\N$
either $\restrict{\CM_{\N}}{A}$ or $\restrict {\CM_\N}{(\N \setminus
A)}$ has the same age as $\CM_\N$. This amounts to saying  that
for each integer $n$, $\restrict{\CM_\N}{[0, n[}$ is
embeddable into either $\restrict{\CM_{\N}}{A}$ or into
$\restrict {\CM_\N}{(\N \setminus A)}$. The embeddability of
$\restrict{\CM_{\N}}{[0, n[}$ into a subset amounts to the
existence of an arithmetic progression $(a+ix)_{i \in n}$ in that
subset.  Van der Waerden's theorem on arithmetic progressions
\cite{graham} ensures the required conclusion.

\noindent Item 2. The additive structure on $V$ is a torsion free 
abelian group and hence satisfies the requirements of
Lemma~\ref{arithmetic}.  Let $A$ be given by Lemma~\ref{arithmetic}
and let $B:= V\setminus A$. The age of $\restrict{\CM_{V}}{A}$
is a proper subset of the age of $\CM_{V}$, because $A$ does
not contain a three element arithmetic progression and hence does not
contain a triple in the relation $\mu_V$. According to
Lemma~\ref{lem:affineembedding}, an embedding $\alpha$ from
$\CM_V$ into itself is an affine map. There does not exist such
an embedding whose range is a subset of $B$ because $B$ contains no
affine line.

\noindent Item 3. Let $V':= \Q^{[R]}$ be the set of functions 
$h:R\rightarrow \Q$ which are $0$ almost everywhere. Let $\delta:
R \rightarrow V'$ be the map defined by $\delta(x)(y):=1$ if $x=y$
and $\delta(x)(y):=0$ otherwise.  $V'$ is a vector space over $\Q$
under the natural addition and scalar multiplication operations, and
let $W$ be the subspace of $V'$ generated by the vectors of the form
$\delta(x)+\delta(z)-2\delta(y)$ such that $(x,y,z)\in
\mu'$. Let $V'/W$ be the quotient of $V'$ by $W$ and let $\rho:
V'\rightarrow V'/W$ be the quotient map.

\begin{claim}\label{claim:embedding}
The map $\rho':=\rho \circ \delta$ is an embedding of $\CR$ into $\CM_{V'/W}$.  
\end{claim}

\noindent {\bf Proof of  Claim~\ref{claim:embedding}.}
We first verify that $(x,y,z)\in \mu'$ if and only if $\blp \rho'(x),
\rho'(y),\rho'(z) \brp \in \mu_{V'/W}$. The ``only if'' part of this equivalence
is immediate: by definition of $W$, $(x,y,z)\in \mu'$ implies
$\delta(x)+\delta(z)-2\delta(y)\in W$.  This amounts to
$\rho \blp \delta(x)+\delta(z)-2\delta(y) \brp =0$, that is
$\rho \blp \delta(x) \brp +\rho \blp \delta(z))-2\rho(\delta(y) \brp =0$ which rewrites as 
$\rho'(x)+\rho'(z)-2\rho'(y)=0$, that is $(\rho'(x), \rho'(y), \rho'(z))\in \mu_{V'/W}$.

\noindent For the ``if'' part, it suffices to show that
$\delta(x)+\delta(z)-2\delta(y)\in W$ implies $(x,y,z)\in
\mu'$. So suppose that $\delta(x)+\delta(z)-2\delta(y)$ is a finite
linear combination $\sum_{i \in n}\lambda_i \blp delta(x_i) +
\delta(z_i)-2\delta(y_i) \brp $ where $(x_i, y_i,z_i)\in \mu'$ 
and $\lambda_i \in \Q$ for each $i \in n$. Let $F:= \{x, y, z, x_i, y_i,
z_i: i \in n\}$, and by hypothesis let $f$ be an isomorphism of
$\CR_{\restriction F}$ into $\CV$. As a map defined on a subset of $R$,
$f$ extends to a linear map $\overline f$ from $V'$ to $V$. As
such it satisfies:
\[ \overline f \blp \delta (x)  \brp +\overline f \blp \delta(z) \brp -2\overline f \blp \delta(y)\brp 
	=\sum_{i \in n}\lambda_i [ \overline f \blp \delta(x_i) \brp +\overline f \blp \delta(z_i) \brp -2\overline f \blp \delta(y_i) \brp ].
\]
Since $f$ preserves $\mu'$, $\blp f(x_i), f(y_i), f(z_i) \brp \in \mu_V$ for
all $i \in n$, hence $\overline f \blp \delta(x_i) \brp +\overline
f \blp \delta(z_i) \brp -2\overline f \blp \delta(y_i) \brp =0$. This yields $\overline
f \blp \delta (x) \brp +\overline f \blp \delta(z) \brp -2\overline f \blp \delta(y) \brp =0$, that
is $\blp f(x), f(y), f(z) \brp \in \mu_V$ from which it follows that
$(x,y,z)\in \mu'$. 

To conclude, it suffices to prove that $\rho'$ is one to one. Let $a,
a'\in R$ such that $\rho'(a)=\rho'(a')$. This means that
$\delta(a)-\delta(a')$ is a finite linear combination
$\sum_{i \in n}\lambda_i \blp \delta(x_i)+\delta(z_i)-2\delta(y_i) \brp $ where
$(x_i, y_i,z_i)\in \mu'$ and $\lambda_i \in \Q$ for each $i \in n$. In
order to prove that this linear combination is zero, we use the same
technique as  above.  Let $F':= \{a, a', x_i, y_i, z_i: i \in n\}$, and
by hypothesis let $f$ be an isomorphism of $\CR_{\restriction F'}$ into
$\CV$. The map $f$ extends to a linear map $\overline f$ from $V'$ to
$V$. It satisfies:
\[
\overline f \blp \delta (a)\brp -\overline f\blp \delta(a') \brp
	=\sum_{i \in n}\lambda_i [ \overline f \blp \delta(x_i) \brp +\overline f \blp \delta(z_i) \brp -2\overline f \blp \delta(y_i) \brp ].
\]
Since $f$ preserves $\mu'$, $\blp f(x_i), f(y_i), f(z_i) \brp \in \mu_V$ for
all $i \in n$, hence $\overline f \blp \delta(x_i) \brp +\overline
f \blp \delta(z_i) \brp -2\overline f\blp \delta(y_i) \brp =0$. Hence $\overline
f \blp \delta(a) \brp = \overline f \blp \delta(a')\brp $  from which it follows that
$a=a'$. This proves our claim.\endproof

Since $R$ is countable, $V'/W$ is countable and hence it is
embeddable into $V$ by some linear map. It follows that
$\CM_{V'/W}$ is embeddable into $\CM_{V}$. Using
Claim~\ref{claim:embedding}, it follows that $\CR$ is embeddable into
$\CM_{V}$. This completes the proof of Theorem \ref{theorem:notweakly}.

\end{proof}

\section{Conclusion}

We have seen that a countable dimensional vector space over an
infinite field is age indivisible, not weakly indivisible.  Although
homogeneous as a relational structure, it has infinitely many
relations. The midpoint structure above is also
age indivisible, not weakly indivisible.  In this case, as a
relational structure, it has a single ternary relation, but is not
homogeneous. It can be made homogeneous by adding relations, but as
above an infinite number is required.

It is natural to impose finiteness conditions in asking for a
countable age indivisible, not weakly indivisible, homogeneous
relational structure. Beside a finite number of relations, one may ask
for an oligomorphic automorphism group (finite number of orbits on
$n$-tuples for each $n$, \cite{Camer}), or an automorphism group of
finite arity (types of $n$-tuples determined by their $r$-tuples for a
fixed $r$, see \cite{Cherl}). We therefore ask:

\begin{question}
Can one impose ``finiteness'' conditions for an age indivisible, not
weakly indivisible, countable homogeneous relational structure?
\end{question}

\end{document}